\documentclass[11pt]{article}

\usepackage[T1]{fontenc}
\usepackage[utf8]{inputenc}
\usepackage{authblk}

\usepackage{amssymb,amsmath,amsthm}
\usepackage[colorlinks=true,linktoc=page,urlcolor=blue,
citecolor=red,linkcolor=blue,pdfpagelabels,bookmarksnumbered,bookmarksopen]{hyperref}
\usepackage{graphicx}
\usepackage[english]{babel}
\usepackage{epsfig}
\usepackage{graphics}

\usepackage{amsthm}
\theoremstyle{plain}

\newtheorem{lemma}{Lemma}[section]

\theoremstyle{definition}

\theoremstyle{plain}
\newtheorem*{conj*}{Conjecture}
\newtheorem*{corollary*}{Corollary}

\newtheorem{prop}[lemma]{Proposition}
\makeatletter
\renewcommand*\l@section[2]{%
  \ifnum \c@tocdepth >\z@
    \setlength\@tempdima{1.5em}%
    \begingroup
      \parindent \z@ \rightskip \@tocrmarg
      \parfillskip -\rightskip
      \leavevmode
      \hangindent\@tempdima \hangafter\@ne
      #1\nobreak\hfil \nobreak\hb@xt@\@pnumwidth{\hss #2}\par
    \endgroup
  \fi}
\makeatother
\theoremstyle{plain}
\newtheorem{thmletter}{Theorem}

\def\div{\mathop{\mathrm{div}}}

\newcommand{\N}{{\mathbb N}}
\newcommand{\R}{{\mathbb R}}

\newcommand{\eps}{\varepsilon}

\newcommand{\DR}{{\cal D}}

\def\bean#1\eean{\begin{eqnarray*}#1\end{eqnarray*}}

\begin{document}

\title{\bfseries\Large Unique continuation, nonexistence and bubbling\\ for the critical $p$-Laplace equation II:\\ the degenerate elliptic case and higher dimensions}
\author{Carlo Mercuri}
\affil{

Dipartimento di Scienze Fisiche, Informatiche e Matematiche\\\small
Universit\`a di Modena e Reggio Emilia\\\small
Via Campi 213/b, 41125 Modena, Italy\\\small
\textit{carlo.mercuri@unimore.it}
}
\date{}

\maketitle
{\renewcommand{\thefootnote}{}}

\begin{abstract}
For every $p>1$, we prove that a nontrivial solution of
\begin{equation*}
-\Delta_pu=f(x,u,\nabla u)
\end{equation*}
in a planar domain, with $|f(x,s,\xi)|\leq C(|s|^{p-1}+|\xi|^{p-1})$, vanishes neither on an open set
nor to infinite order at a point. This allows us to prove unique continuation properties of axially
symmetric solutions to $p$-Laplace equations in higher dimensions $N\geq3$, with axially symmetric
$p$-harmonic functions as a particular case. Unlike in our previous work \cite{M}, the analysis is
based here on the regularity and quasiconformal properties of the whole stress field
$|\nabla u|^{p-2}\nabla u$ in the plane, taking into account the more delicate regularity theory for
$p>2$. This allows us to prove Liouville-type and bubble profile decomposition results for critical
problems in dimensions $N\geq3$ and in the whole range $1<p<N$, under some symmetry conditions.
\end{abstract}

\setcounter{tocdepth}{1}
\tableofcontents

\section[Introduction: a quasiconformal reduction principle and main results]{Introduction: a quasiconformal reduction principle\\ and main results}
\label{sec:intro}

In this paper we complete the study started in \cite{M} of unique continuation properties for
$p$-Laplace equations, considering the more general structure
\begin{equation}
\label{eq:eq-f}
-\Delta_pu=f(x,u,\nabla u)\qquad\mbox{in }\Omega ,
\end{equation}
where $\Omega\subset\R^N$ is a domain, allowing the source term to depend also on the gradient. In particular
we consider every
$N\geq2$ and $p>1$.

In \cite{M} we showed that solutions of
\begin{equation*}
-\Delta_pu=f(u)
\end{equation*}
in planar domains, for $1<p<2$ and $f$ continuous with
\begin{equation*}
|f(s)|\leq C_M|s|^{p-1}\qquad\mbox{for }|s|\leq M ,
\end{equation*}
which vanish on an open set or to infinite order at a point vanish identically, and we derived from
this nonexistence and global compactness results for the critical $p$-Laplace equation in the plane.

As anticipated in \cite{M}, and as we elaborate on here, equations of the
type \eqref{eq:eq-f} occur, for instance, when
considering axially symmetric solutions in dimension $N\geq3$. In the axial coordinates $(r,x_N)$,
$r=|x'|$, a solution $u(x)=U(r,x_N)$ of
\begin{equation*}
-\Delta_pu=f(x,u)
\end{equation*}
in a domain of $\R^N$ is reduced, on the half-plane $\{r>0\}$, to
\begin{equation*}
-\Delta_pU=f(r,x_N,U)+\frac{N-2}r\,|\nabla U|^{p-2}\partial_rU .
\end{equation*}
Moreover, critical equations in dimension
$N\geq3$ allow $p$ to vary naturally from the singular $p<2$ to the degenerate elliptic $p>2$ regimes. It is therefore desirable to extend the unique continuation results
of \cite{M} to $p>2$ and to source terms depending on the gradient. To this aim the mathematical challenge to be overcome is mainly twofold: the regularity of the
solutions, which is more delicate for $p>2$, and the method of \cite{M}, which is robust enough to
be adapted to this more general setting, yet subtle enough to require a more careful analysis of the
reduction to a Beltrami-type system, carried out here for the whole stress field, as we explain below,
outlining a general method.

For every $p>1$, we prove in Theorem \ref{thm:plane} that a nontrivial solution of \eqref{eq:eq-f}
in a planar domain, with
\begin{equation*}
|f(x,s,\xi)|\leq C_M\big(|s|^{p-1}+|\xi|^{p-1}\big)\qquad\mbox{for }|s|\leq M,\ |\xi|\leq M ,
\end{equation*}
vanishes neither on an open set nor to infinite order at a point. For
$p$-harmonic functions in the plane unique continuation holds
by the celebrated work of Manfredi \cite[Theorem 1 and Corollary 1]{Manfredi}, as the complex gradient
$h=\partial_zu$ is quasiregular.

Guo and Kar \cite{GuoKar} have considered strong unique continuation and quantitative estimates of
the order of vanishing and of the rate of decay at infinity for
\begin{equation*}
-\Delta_pu=A\cdot|\nabla u|^{p-2}\nabla u
\end{equation*}
in planar domains, with $A$ locally Lipschitz. Our results therefore include strong unique
continuation in this case, also in the presence of a zeroth-order term. As we have already observed
in \cite{M}, away from the axis the equation of an axially symmetric $p$-harmonic
function, and in general of an axially symmetric solution of \eqref{eq:eq-f}, in dimension
$N\geq3$, reduces to an equation of the type \eqref{eq:eq-f} in the plane. Unique continuation for these
functions is obtained in Theorem \ref{thm:axial} as a particular case, together with the isolatedness of
their critical orbits on the nodal set, each of index $-n$.

Unlike in \cite{M}, the analysis is based here on the regularity and quasiconformal properties of the
whole stress field $|\nabla u|^{p-2}\nabla u$ in the plane. In \cite{M}, for $1<p<2$, the first-order
system for the complex gradient $h$ is derived from the non-divergence form of the equation
\begin{equation*}
a_{ij}\,\partial_{ij}u=-f(x,u,\nabla u)\,|\nabla u|^{2-p},
\end{equation*}
with
\begin{equation*}
a_{ij}=\delta_{ij}+(p-2)\,e_ie_j,\qquad e=\frac{\nabla u}{|\nabla u|},
\end{equation*}
away from the critical set, and $a_{ij}=\delta_{ij}$ on it; the zeroth-order term produces in the
system a term $\beta u$, not involving the unknown $h$, which is removed by the transform
\begin{equation*}
\mathcal R=h-vu ,
\end{equation*}
with $v$ a solution to a quadratic Riccati-type equation. Two properties of the range $1<p<2$ make this work.
The second derivatives of $u$ are square summable, which gives a meaning to the non-divergence form
almost everywhere, since by a classical observation of Stampacchia
\begin{equation*}
D^2u=0\qquad\mbox{a.e.\ on }\{\nabla u=0\},
\end{equation*}
see \cite[(13)]{M}, and places $h$ in $W^{1,2}_{{\textrm{loc}}}$, as required of a generalised solution
in the sense of Bojarski \cite[Section 2]{Bojarski}. Moreover, the growth condition and Young's inequality
give, locally in $\Omega$,
\begin{equation*}
|a_{ij}\,\partial_{ij}u|\leq C\big(|u|+|\nabla u|\big),
\end{equation*}
which bounds the coefficients of the system for $h$, see \cite[(15)]{M}. For $p>2$ this breaks down:
the right-hand side of the non-divergence form is no longer controlled by $|u|+|\nabla u|$ near the
critical set and, already for $p\geq3$, the second derivatives may fail to be square summable, as pointed out
in Damascelli and Sciunzi \cite[Proposition 2.2]{DamascelliSciunzi}, while the lower-order coefficients of
the system for $h$ may fail to belong to any $L^s_{{\textrm{loc}}}$ with $s>2$, as required by the theory of
Bojarski on which the representation of \cite[Lemma 4.1]{M} is based.

We bypass these issues by writing the equation as a Beltrami-type system for the pair
\begin{equation*}
h_p=|h|^{p-2}h,\qquad u_p=|u|^{p-2}u ,
\end{equation*}
in place of $(h,u)$, where $h_p$ is the complex representation of the full stress field. Our equation \eqref{eq:eq-f}
then reads as
\begin{equation*}
\begin{cases}
\partial_{\bar z}h_p-\mu\,\partial_zh_p-\nu\,\overline{\partial_zh_p}=\alpha h_p+\beta u_p,\\
\partial_zu_p=a\,h_p+b\,u_p ,
\end{cases}
\end{equation*}
with $|\mu|+|\nu|\leq\frac{p-2}p<1$, see Proposition \ref{prop:system}, and its lower-order
coefficients are locally bounded, as required by Bojarski's theory. As we shall see, this is a
consequence of the fact that the pair $(h_p,u_p)$ shares the homogeneity of the operator, of degree
$p-1$: for $p\geq2$, locally in $\Omega$,
\begin{equation*}
|f(x,u,\nabla u)|\leq C\big(|h_p|+|u_p|\big),\qquad |\partial_zu_p|\leq(p-2)|u_p|+|h_p| ,
\end{equation*}
which bound $\alpha,\beta$ and $a,b$, respectively. The first estimate plays the role of the key
estimate \cite[(15)]{M}, while the second relaxes, and plays the same role as, the identity
$\partial_zu=h$, which couples $u$ and $h$ in \cite{M}.
The only regularity used is that the stress field lies in $W^{1,2}_{{\textrm{loc}}}$: this goes back to
Bojarski and Iwaniec for $p$-harmonic functions and $p\geq2$, see e.g.\
Lindqvist \cite[Theorem 4.1]{Lindqvist}, and to Damascelli and
Sciunzi \cite[Corollary 2.2]{DamascelliSciunzi} for positive solutions, and it holds in the
present general setting by Cianchi and Maz'ya \cite[Theorem 2.1]{CianchiMazya}. By the
transform
\begin{equation*}
\mathcal R=h_p-v\,u_p ,
\end{equation*}
for a suitable $v$ as in \cite{M}, we are then able to remove the term $\beta u_p$, since $\mathcal R$
satisfies the homogeneous system
\begin{equation*}
\partial_{\bar z}\mathcal R=\mu\,\partial_z\mathcal R+\nu\,\overline{\partial_z\mathcal R}+A\mathcal R+B\overline{\mathcal R} ,
\end{equation*}
with $A$ and $B$ bounded. The representation theorem of Bers and Nirenberg, in the form of
Bojarski \cite[Theorem 4.4]{Bojarski}, then shows that the zeros of $\mathcal R$ are isolated unless
$\mathcal R$ vanishes identically, see Lemma \ref{lem:repr-W}. Finally, as in \cite{M}, this and a
Caccioppoli-type estimate allow us to prove that a nontrivial solution vanishes neither on an open set
nor to infinite order at a point.

Unique continuation in its weak form allows us to prove Liouville-type and bubble
profile decomposition results for critical problems in the whole range $1<p<N$, extending \cite[Theorems B, C
and D]{M} to axially symmetric solutions in dimension $N\geq3$; these are Theorems \ref{thm:nonex}
and \ref{thm:compact} below.

\subsection*{A quasiconformal reduction principle}

We may call the method we follow in \cite{M} and here a quasiconformal reduction principle: a
Beltrami-type system for a solution of our equation and its gradient is reduced, by a transform, to a
homogeneous one, so that the representation theorem of Bers and Nirenberg holds, in a
generalised form due to Bojarski, yielding the isolatedness of the zeros
and hence a weak unique continuation property. This improves to the strong form of unique continuation
provided a Caccioppoli-type inequality is available, which is typically the case if the equation is in
divergence form.

Our equation is written as the Beltrami-type system
\begin{equation*}
\begin{cases}
\mathcal LW=\alpha W+\beta U,\\
\partial_zU=aW+bU ,
\end{cases}
\qquad \mathcal L=\partial_{\bar z}-\mu\,\partial_z-\nu\,\overline{\partial_z},
\end{equation*}
for a pair $(W,U)$ of a complex and a real unknown, vanishing where $(\nabla u,u)$ does, with
$|\mu|+|\nu|\leq k<1$: $(W,U)=(h,u)$ in \cite{M}, where $1<p<2$, and $(W,U)=(h_p,u_p)$ here,
where $p\geq2$. The growth condition on $f$ bounds $\alpha$ and $\beta$, so that the source term may
depend on the gradient. The divergence form of the equation also allows the system to be derived for
the stress field, with no a priori information on the second derivatives of $u$. In the plane, the
$\partial_{\bar z}$-derivative of a vector field $F=(F_1,F_2)$ can be represented as
\begin{equation*}
\partial_{\bar z}(F_1-iF_2)=\frac12\big(\div F-i\,\mathrm{curl}\,F\big),\qquad \mathrm{curl}\,F=\partial_1F_2-\partial_2F_1 .
\end{equation*}
For the stress field $F=|\nabla u|^{p-2}\nabla u$, represented by $h_p$ in its complex form, equation \eqref{eq:eq-f}
prescribes the divergence,
\begin{equation*}
\div F=-f(x,u,\nabla u) ,
\end{equation*}
while, as $\nabla u$ is a function of $F$ with zero curl in the sense of distributions, the curl of $F$
can be controlled by a fraction of its first derivatives, which produces a Beltrami-type structure, see
Proposition \ref{prop:system}.

Finally, the transform
\begin{equation*}
\mathcal R=W-vU
\end{equation*}
allows us to reduce to a Beltrami-type system without the term $\beta U$, so that the representation
theorem of Bers and Nirenberg can be applied to $\mathcal R$. This is achieved by solving a quadratic
Riccati-type equation for $v$ with a classical Banach--Caccioppoli contraction mapping argument.

\subsection*{Unique continuation in the plane for all $p>1$}

Let $\Omega\subset\R^N$ be open, $N\geq2$, and $p>1$, and let $f=f(x,s,\xi)$ be measurable in $x$ and
continuous in $(s,\xi)$. We say that $u\in C^{1,\alpha}_{{\textrm{loc}}}(\Omega)$ is a weak solution
of \eqref{eq:eq-f} if $f(x,u,\nabla u)\in L^1_{{\textrm{loc}}}(\Omega)$ and
\begin{equation*}
\int_\Omega|\nabla u|^{p-2}\nabla u\cdot\nabla\varphi\,dx=\int_\Omega f(x,u,\nabla u)\,\varphi\,dx
\qquad\mbox{for every }\varphi\in C^\infty_c(\Omega) .
\end{equation*}
We set
\begin{equation*}
Z=\{\nabla u=0\},\qquad \Sigma_0=Z\cap\{u=0\} .
\end{equation*}
For a closed path $\gamma$ in the plane and $z$ not lying on $\gamma$, the index of $z$ with respect to
$\gamma$ is the integer
\begin{equation*}
\mathrm{Ind}_\gamma(z)=\frac1{2\pi i}\int_\gamma\frac{dw}{w-z} ,
\end{equation*}
see e.g.\ Rudin \cite{Rudin}; the definition extends, by
uniform approximation, to closed curves. If a continuous vector field $X=(X_1,X_2)$, which we
identify with $X_1+iX_2$, has an isolated zero at $z_0$, its index at $z_0$ is
\begin{equation*}
\mathrm{ind}_{z_0}X=\mathrm{Ind}_{X\circ\gamma_r}(0),\qquad \gamma_r(t)=z_0+re^{it},\quad 0\leq t\leq2\pi ,
\end{equation*}
where $\gamma_r$ describes $\partial B_r(z_0)$ counterclockwise, for $r>0$ small; this does not depend on
$r$. We use it for $X=\nabla u$ at the points of $\Sigma_0$:
by Theorem \ref{thm:plane}(iii) below, each of them is an isolated point of the whole critical set $Z$,
that is, no other critical point of $u$, at any level, accumulates at it, so that the index is
defined, and it equals $-n$ for some integer $n\geq1$. It is easy to construct examples, even for $p=2$, in which isolatedness fails at critical
points at nonzero levels, that is at points of $Z\setminus\Sigma_0$.

\begin{thmletter}[Unique continuation in the plane]
\label{thm:plane}
Let $p>1$, let $\Omega\subset\R^2$ be open and connected, let $f=f(x,s,\xi)$ be measurable in $x$,
continuous in $(s,\xi)$, and such that, for every $M>0$,
\begin{equation}
\label{eq:f-growth}
|f(x,s,\xi)|\leq C_M\big(|s|^{p-1}+|\xi|^{p-1}\big)\qquad\mbox{for a.e.\ }x\in\Omega,\ |s|\leq M,\ |\xi|\leq M ,
\end{equation}
and let $u\in C^{1,\alpha}_{{\textrm{loc}}}(\Omega)$ be a weak solution of \eqref{eq:eq-f}.
\begin{enumerate}
\item[(i)] If $u$ vanishes on a nonempty open subset of $\Omega$, then $u\equiv0$ in $\Omega$.
\item[(ii)] If $u\not\equiv0$, then for every $z_0\in\Omega$ there are $m\geq0$ and $c,\rho_0>0$ with
\begin{equation}
\label{eq:finite-order}
\int_{B_\rho(z_0)}|u|^p\,dx\ \geq\ c\,\rho^{m}\qquad\mbox{for }0<\rho<\rho_0 ,
\end{equation}
so that $u$ vanishes to a finite order at each of its zeros. In particular, a solution which vanishes
to infinite order at one point, in the sense that the integral in \eqref{eq:finite-order} is
$O(\rho^j)$ as $\rho\to0$ for every $j\in\mathbb N$, vanishes identically.
\item[(iii)] If $u\not\equiv0$, then every point $z_0$ of $\Sigma_0=Z\cap\{u=0\}$ is an isolated point
of $Z$, and the index of $\nabla u$ at $z_0$ is $-n$, for some integer $n\geq1$.
\end{enumerate}
\end{thmletter}

The statements regarding weak and strong unique continuation
extend \cite[Theorem A]{M}, proved there for $1<p<2$ and $f=f(u)$. In the singular regime $1<p<2$ they
are, in essence, those of \cite{M}, whose proof accommodates the general right-hand sides
of \eqref{eq:eq-f}, as observed in \cite[Remark 1.1]{M}; the degenerate regime $p\geq2$ is based on
handling the whole stress field, and both are proved in Section \ref{sec:uc}. Assertion (iii) extends \cite[Corollary 4.3]{M}, where $\Sigma_0$ is shown to
consist of isolated points for $1<p<2$ and $f=f(u)$, and answers for every $p>1$ the question raised
in \cite[Related question 2)]{M} on the index of $\nabla u$ at these points; this is proved in
Section \ref{sec:index}.

For $p$-harmonic functions in the plane the critical points of a nonconstant solution are
isolated \cite[Corollary 1]{Manfredi}, and near each of them it holds that
\begin{equation*}
h=\chi^n ,
\end{equation*}
where $\chi$ is a quasiconformal homeomorphism vanishing there and $n\geq1$ is an integer, so that the index of
$\nabla u$ is $-n$; see Iwaniec and Manfredi \cite[(12)--(15)]{IwaniecManfredi}. For $p=2$ and
$f=f(u)$ the formula is in Alessandrini and Magnanini \cite[Theorem 3.5]{AlessandriniMagnanini},
under assumptions which include the isolatedness of the critical points. Here only the modified stress
field $\mathcal R$ has, near each point of $\Sigma_0$, the form $\chi^n$ up to a nonvanishing factor, and we
transfer the isolatedness and the index from $\mathcal R$ to $h$ using, as in the proof
of \cite[Proposition 3.3]{GuoKar}, the fact that quasiconformal homeomorphisms map small circles onto
curves lying between two concentric circles of comparable radii, see Section \ref{sec:index}.

In particular, for every $p>1$, a solution
$u\not\equiv0$ of
\begin{equation*}
-\Delta_pu=\lambda|u|^{p-2}u
\end{equation*}
on a planar domain, such as an eigenfunction of the $p$-Laplacian, has a nodal set with empty
interior, has a finite order of vanishing at each of its zeros, and its critical points on the nodal set
are isolated in the critical set, each of index $-n$ for some integer $n\geq1$. On the eigenvalue
problem for the $p$-Laplacian we refer to Lindqvist \cite{LindqvistEig}.

\subsection*{Axial symmetry and higher dimensions}

Symmetry allows the planar results to be lifted to higher dimensions: when a compact group of
isometries of $\R^N$ has a two-dimensional orbit space, the equation for invariant solutions becomes,
away from the singular orbits, such as the points of the axis for the rotations about it, a planar
equation with first-order terms of the type considered above.

We deal here with the rotations about an axis. Let $N\geq3$, $x=(x',x_N)\in\R^{N-1}\times\R$,
$r=|x'|$, $A=\{x'=0\}$ the $x_N$-axis, and $G=SO(N-1)\times\{1\}$ the group of rotations about $A$.
Recall that the $G$-orbit of $x$ is the set $Gx=\{gx:\ g\in G\}$. A set $\Omega$ is axially symmetric if
$Gx\subset\Omega$ for every $x\in\Omega$, and a function is axially symmetric if it is constant on every
$G$-orbit, that is, if
\begin{equation*}
u(x)=U(r,x_N) .
\end{equation*}
For $x_0\notin A$ the orbit $Gx_0$ is a sphere of dimension $N-2$ and, since
$\operatorname{dist}(x,Gx_0)=\big|(|x'|,x_N)-(|x_0'|,x_{0,N})\big|$, its tubular neighbourhoods
\begin{equation*}
B_\rho(Gx_0)=\big\{x:\ \operatorname{dist}(x,Gx_0)<\rho\big\}
\end{equation*}
correspond to the discs of radius $\rho$ about $(|x_0'|,x_{0,N})$ in the half-plane $\{r>0\}$. If
$u\in C^1$ is axially symmetric and $x_0\in Z\setminus A$, we
say that $Gx_0$ is an \emph{isolated orbit} of $Z$ when $B_\rho(Gx_0)\cap Z=Gx_0$ for some $\rho>0$,
that is, when $(|x_0'|,x_{0,N})$ is an isolated critical point of $U$, and we call the index of $\nabla U$
at that point the index of the orbit.

\begin{thmletter}[Unique continuation under axial symmetry]
\label{thm:axial}
Let $N\geq3$, $p>1$, let $\Omega\subset\R^N$ be an axially symmetric domain, let $f$ satisfy
\eqref{eq:f-growth}, and let $u\in C^{1,\alpha}_{{\textrm{loc}}}(\Omega)$ be an axially symmetric weak
solution of \eqref{eq:eq-f}.
\begin{enumerate}
\item[(i)] If $u$ vanishes on a nonempty open subset of $\Omega$, then $u\equiv0$ in $\Omega$.
\item[(ii)] If $u\not\equiv0$, then for every orbit $Gx_0\subset\Omega\setminus A$ there are $m\geq0$ and
$c,\rho_0>0$ with
\begin{equation*}
\int_{B_\rho(Gx_0)}|u|^p\,dx\ \geq\ c\,\rho^{m}\qquad\mbox{for }0<\rho<\rho_0 ,
\end{equation*}
so that $u$ vanishes to a finite order on each orbit of zeros away from the axis.
\item[(iii)] If $u\not\equiv0$, then $\Sigma_0\setminus A$ is a union of isolated orbits of $Z$, each
of index $-n$ for some integer $n\geq1$.
\end{enumerate}
\end{thmletter}

Theorem \ref{thm:axial} applies to axially symmetric $p$-harmonic functions, $f=0$. For instance, an
axially symmetric $p$-harmonic function in $\R^3$ which vanishes in the half-space $\{x_3<0\}$ vanishes
everywhere. This gives, under symmetry, a partial answer to a well-known open problem, see
e.g.\ Lindqvist \cite[Chapter 11]{Lindqvist}. Theorem \ref{thm:axial} is proved in
Section \ref{sec:axial}, where the profile $U$ is shown to solve, in the half-plane $\{r>0\}$, a planar
equation to which Theorem \ref{thm:plane} applies.

This low-dimensional reduction in the presence of axial symmetry is reminiscent of the work of Bers and Gelbart \cite[Section 5]{BersGelbart} in linear potential
theory, where axially symmetric potentials are described by systems of Cauchy--Riemann type in the
meridian half-plane. For compact groups of isometries of $\R^N$ having a two-dimensional orbit space,
with calculations similar to those of Section \ref{sec:axial}, it is possible to show that, away from
the singular orbits, invariant solutions satisfy a reduced equation of the type considered in
Theorem \ref{thm:plane}.

\subsection*{Nonexistence for the critical equation}

Let now $1<p<N$, $p^*=\frac{Np}{N-p}$, $H=\{x\in\R^N:\ x_N>0\}$, and let $\DR_0^{1,p}(H)$ and the weak
solutions of $-\Delta_pu=|u|^{p^*-2}u$ in $H$, or in a bounded domain $\Omega$ with $W^{1,p}_0(\Omega)$
in place of $\DR_0^{1,p}(H)$, be as in \cite[Section 1.2]{M}. The axis $A$ is normal to
$\partial H$, so that $H$ is axially symmetric.

\begin{thmletter}[Nonexistence under axial symmetry]
\label{thm:nonex}
Let $N\geq3$ and $1<p<N$.
\begin{enumerate}
\item[(i)] Every axially symmetric weak solution $u\in\DR_0^{1,p}(H)$ of
\begin{equation*}
-\Delta_pu=|u|^{p^*-2}u\qquad\mbox{in }H
\end{equation*}
vanishes identically.
\item[(ii)] Let $\Omega\subset\R^N$ be a bounded axially symmetric domain with $C^2$ boundary,
star-shaped with respect to a point $x_0\in\Omega$:
\begin{equation*}
(x-x_0)\cdot\nu(x)\geq0\qquad\mbox{for every }x\in\partial\Omega ,\end{equation*}
where $\nu$ is the outer unit normal. Then every axially symmetric weak solution $u\in W^{1,p}_0(\Omega)$ of
\begin{equation*}
-\Delta_pu=|u|^{p^*-2}u\qquad\mbox{in }\Omega
\end{equation*}
vanishes identically.
\end{enumerate}
\end{thmletter}

Part (i) removes, under axial symmetry, the sign condition in the nonexistence theorem of
the author and Willem \cite[Theorem 1.1]{MercuriWillem}. Without the symmetry assumption, it would
extend to the $p$-Laplacian, in the whole range $1<p<N$, the theorem of Esteban and
Lions \cite{EstebanLions} on the half-space, and remove the sign condition from the representation
theorem of the author and Willem \cite[Theorem 1.2]{MercuriWillem}, completing the generalisation of
Struwe's theorem \cite{Struwe}. The latter is thus reduced to a symmetry question.

Part (ii), in the spirit of Pohozaev, extends to axially symmetric solutions in star-shaped domains
the nonexistence of radial solutions in a ball of \cite[Proposition 1.2]{FarinaMercuriWillem}.
Theorem \ref{thm:nonex} is proved in Section \ref{sec:nonex}.

\subsection*{Global compactness}

Let $\Omega\subset\R^N$ be a smooth bounded axially symmetric domain, and let $a\in L^{N/p}(\Omega)$
be axially symmetric. We regard $W^{1,p}_0(\Omega)$ as a subspace of $\DR^{1,p}(\R^N)$, by extension by
zero, with the norm $\|u\|=\|\nabla u\|_{L^p}$, we denote by $W^{1,p}_{0,G}(\Omega)$ the subspace of its
axially symmetric functions, and we set
\begin{align*}
J_\Omega(u)&=\int_\Omega\frac{|\nabla u|^p}p+a\frac{|u|^p}p-\frac{|u|^{p^*}}{p^*}\,dx,\\
J_{\R^N}(u)&=\int_{\R^N}\frac{|\nabla u|^p}p-\frac{|u|^{p^*}}{p^*}\,dx .
\end{align*}
The following result states that, under axial
symmetry, Palais--Smale sequences lose compactness only through bubbles which concentrate at points of
the axis.

\begin{thmletter}[Global compactness under axial symmetry]
\label{thm:compact}
Let $N\geq3$, $1<p<N$, let $\Omega$ and $a$ be as above, and let
$(u_n)_{n\in\N}\subset W^{1,p}_{0,G}(\Omega)$ be a Palais--Smale sequence for $J_\Omega$ at a level $c$.
Then, passing if necessary to a subsequence, there are an integer $k\geq0$, an axially symmetric solution
$v_0\in W^{1,p}_0(\Omega)$ of
\begin{equation*}
-\Delta_pu+a|u|^{p-2}u=|u|^{p^*-2}u\qquad\mbox{in }\Omega ,
\end{equation*}
nontrivial axially symmetric solutions $v_1,\dots,v_k\in\DR^{1,p}(\R^N)$ of
\begin{equation*}
-\Delta_pu=|u|^{p^*-2}u\qquad\mbox{in }\R^N ,
\end{equation*}
and sequences $(\eps_n^i)_{n\in\N}\subset(0,\infty)$, $(y_n^i)_{n\in\N}\subset\Omega\cap A$,
$i=1,\dots,k$, with
\begin{equation*}
\eps_n^i\to0,\qquad\frac{\operatorname{dist}(y_n^i,\partial\Omega)}{\eps_n^i}\to\infty ,
\end{equation*}
such that, as $n\to\infty$,
\begin{gather*}
\Big\|u_n-v_0-\sum_{i=1}^k(\eps_n^i)^{\frac{p-N}p}\,v_i\Big(\frac{\cdot-y_n^i}{\eps_n^i}\Big)\Big\|\to0,\\
\|u_n\|^p\to\|v_0\|^p+\sum_{i=1}^k\|v_i\|^p,\\
J_\Omega(v_0)+\sum_{i=1}^kJ_{\R^N}(v_i)=c .
\end{gather*}
In particular, the concentration may occur only on the axis, with no limiting profile-solution to the critical
half-space equation occurring at all. Moreover, if $\Omega$ does not meet the axis, then $k=0$ and $u_n\to v_0$
strongly in $W^{1,p}_0(\Omega)$. \end{thmletter}

Since $G$ acts on $W^{1,p}_0(\Omega)$ by isometries and $J_\Omega$ is invariant, its derivative
$J_\Omega'$ is equivariant, see e.g.\ Willem \cite[Definition 1.27]{Willem}.
Hence, as in the principle of symmetric criticality of Palais \cite{Palais}, the Palais--Smale
sequences of the restriction of $J_\Omega$ to $W^{1,p}_{0,G}(\Omega)$ are Palais--Smale sequences of
$J_\Omega$, to which Theorem \ref{thm:compact} applies. Theorem \ref{thm:compact} is proved in
Section \ref{sec:compact}.

Without symmetry, the analogue of Theorem \ref{thm:compact} is \cite[Theorem D]{M} for $N=2$
and $1<p<2$, and \cite[Theorem 1.2]{MercuriWillem} for $1<p<N$, under the hypothesis that the
negative parts of the $u_n$ tend to zero in $L^{p^*}(\Omega)$, while for $p=2$ and $N\geq3$ it is the
classical global compactness theorem of Struwe \cite{Struwe}.

From Struwe's theorem, Clapp and Pacella \cite[Lemma 4.1(i)]{ClappPacella} obtain the Palais--Smale condition for sequences invariant
under a closed subgroup $\Gamma$ of $O(N)$ below the level
\begin{equation*}
\frac\ell N\,S^{N/2},
\end{equation*}
where $S$ is the best Sobolev constant and $\ell$ is the minimal cardinality of the $\Gamma$-orbits in
$\R^N\setminus\{0\}$.

Half-space profiles are ruled out by a Liouville-type theorem on the half-space, so
that, by Theorem \ref{thm:nonex}, only limiting profiles in $\R^N$ occur, even when the boundary of the
domain meets the axis, as the boundary of a ball does at the poles. As in \cite[Theorem E]{M}, arguing as in
Section \ref{sec:compact}, we observe by Theorem \ref{thm:nonex}(i) that the limiting half-space profiles
in the decomposition of Chernysh \cite{Chernysh} for weighted critical equations do not occur for axially
symmetric Palais--Smale sequences.

When the domain does not meet the axis, as for a solid torus, the Palais--Smale condition holds at
every level, as observed in the celebrated works of P.-L. Lions \cite{LionsSym}, Clapp \cite{Clapp},
and Hebey and Vaugon \cite[Theorem 9.1]{Hebey}, as a consequence of the compactness of the embedding
of $W^{1,p}_{0,G}(\Omega)$ into $L^{p^*}(\Omega)$. In the representation
theorems \cite[Theorem 2.2]{MercuriPacella} and \cite[Theorem 2.5]{ClappLopezRios} for sign-changing
Palais--Smale sequences, limiting profiles living on half-spaces cannot be
excluded, see \cite[Remark 2.5]{MercuriPacella}, since unique continuation is not known for $p\neq2$ in dimension $N\geq3$;
under axial symmetry, they are however ruled out by Theorem \ref{thm:nonex}.

\section{A Beltrami-type system for the complex stress field}
\label{sec:stress}

In this section we show that, in the plane, equation \eqref{eq:eq-f} can be reduced to a homogeneous
Beltrami-type system for a perturbation of the complex stress field,
\begin{equation*}
\mathcal R=h_p-v\,u_p ,
\end{equation*}
which vanishes on $\Sigma_0$ and whose zeros, by the representation theorem of Bers and Nirenberg, are
isolated away from the interior of the nodal set of $u$. We start with the regularity properties of the stress field, for all
$N\geq2$, which are needed throughout the paper.

\subsection*{Regularity of the stress field}
\label{sec:reg}

Let $\Omega\subset\R^N$ be open, $N\geq2$, let $f$ satisfy \eqref{eq:f-growth}, let $u$ be a solution
of \eqref{eq:eq-f}, and set
\begin{equation*}
F=|\nabla u|^{p-2}\nabla u,\qquad g=f(x,u,\nabla u) .
\end{equation*}
The regularity needed for our reduction to a Beltrami-type system follows from a result of Cianchi
and Maz'ya \cite[Theorem 2.1]{CianchiMazya}.

\begin{lemma}[Regularity of the stress field]
\label{lem:stress}
Let $p>1$ and let $u$ be a solution of \eqref{eq:eq-f}. Then
\begin{equation}
\label{eq:stress-reg}
F\in W^{1,2}_{{\textrm{loc}}}(\Omega),\qquad \div F=-g\quad\mbox{a.e.\ in }\Omega ,
\end{equation}
and
\begin{equation}
\label{eq:on-Z}
\partial_iF=0\quad(i=1,\dots,N)\qquad\mbox{and}\qquad g=0\qquad\mbox{a.e.\ on }Z .
\end{equation}
\end{lemma}

\begin{proof}
Since the theorem of Cianchi and Maz'ya applies to the $p$-Laplacian, with $a(t)=t^{p-2}$ in their
notation, we have $F\in W^{1,2}_{{\textrm{loc}}}(B)$ on every ball $B\subset\subset\Omega$. In fact,
by their result, it suffices to find $g_k\in C^\infty_c(B)$, converging to $g$ in some $L^q(B)$, and weak
solutions $u_k$ of $-\Delta_pu_k=g_k$ in $B$, such that
\begin{gather*}
u_k\to u\quad\mbox{and}\quad\nabla u_k\to\nabla u\quad\mbox{a.e.\ in }B,\tag{a}\\
\int_{B'}|\nabla u_k|^{p-1}dx\to\int_{B'}|\nabla u|^{p-1}dx\quad\mbox{for every ball }B'\subset\subset B ,\tag{b}
\end{gather*}
see \cite[(5.1)--(5.2)]{CianchiMazya}. Let $g_k\in C^\infty_c(B)$ converge to $g$, which is bounded on
$B$, in $L^q(B)$ for every $q<\infty$, and let $u_k$ be the weak solutions of
\begin{equation*}
-\Delta_pu_k=g_k\quad\mbox{in }B,\qquad u_k-u\in W^{1,p}_0(B) ;
\end{equation*}
we show that these conditions hold along a subsequence.

Testing with $u_k-u$, we see that $u_k$ is bounded in $W^{1,p}(B)$, and
\begin{equation*}
\int_B\big(|\nabla u_k|^{p-2}\nabla u_k-|\nabla u|^{p-2}\nabla u\big)\cdot\nabla(u_k-u)\,dx
=\int_B(g_k-g)(u_k-u)\,dx\ \to\ 0 .
\end{equation*}
By the monotonicity properties of the $p$-Laplacian, see e.g.\ Lindqvist \cite[Chapter 12]{Lindqvist},
it is now standard to conclude that $\nabla u_k\to\nabla u$ in $L^p(B)$, hence
$u_k\to u$ in $W^{1,p}(B)$. Up to a subsequence, $u_k\to u$ and $\nabla u_k\to\nabla u$ a.e.\ in $B$,
and $|\nabla u_k|^{p-1}\to|\nabla u|^{p-1}$ in $L^{p/(p-1)}(B)$, so that (a) and (b) hold. So, we may conclude that $F\in W^{1,2}_{{\textrm{loc}}}(\Omega)$, and the weak form of the
equation holds almost everywhere, namely \eqref{eq:stress-reg}.

By Stampacchia's classical result, see \cite[Lemma 7.7]{GilbargTrudinger}, we have $\partial_iF=0$
a.e.\ on $Z$ for every $i$, hence $g=-\div F=0$ a.e.\ on $Z$. This concludes the proof.
\end{proof}

\subsection*{The Beltrami-type system}

Throughout the rest of the section, $N=2$ and $p\geq2$. We write $z=x+iy$ and use the complex derivatives of
the Cauchy--Riemann equations,
\begin{equation*}
\partial_z=\frac12\big(\partial_x-i\,\partial_y\big),\qquad \partial_{\bar z}=\frac12\big(\partial_x+i\,\partial_y\big),
\end{equation*}
so that a function $w$ is holomorphic exactly when $\partial_{\bar z}w=0$. The complex gradient of $u$ is
$h=\partial_zu=\frac12(u_x-iu_y)$, and we set
\begin{equation}
\label{eq:pair}
\begin{aligned}
h_p&=|h|^{p-2}h=2^{1-p}(F_1-iF_2), &\qquad u_p&=|u|^{p-2}u ,\\
|h_p|&=2^{1-p}|\nabla u|^{p-1}, &\qquad |u_p|&=|u|^{p-1} .
\end{aligned}
\end{equation}
To state and prove the next proposition, we define the coefficients of the system for the pair
$(h_p,u_p)$. Away from $\Sigma_0$,
\begin{equation}
\label{eq:ab-alphabeta}
(a,b)=\frac{\partial_zu_p}{|h_p|^2+u_p^2}\,\big(\overline{h_p},u_p\big),\qquad
(\alpha,\beta)=-\frac{g}{2^p(|h_p|^2+u_p^2)}\,\big(\overline{h_p},u_p\big),
\end{equation}
and away from $Z$,
\begin{equation}
\label{eq:munu}
\mu=-\frac{p-2}{2p}\,\frac{\overline{h_p}}{h_p},\qquad \nu=\frac{p-2}{2p}\,\frac{h_p}{\overline{h_p}} ;
\end{equation}
we set $a=b=\alpha=\beta=0$ on $\Sigma_0$, and $\mu=\nu=0$ on $Z$.

\begin{prop}[The Beltrami-type system]
\label{prop:system}
Let $N=2$, $p\geq2$, let $u$ be a solution of \eqref{eq:eq-f}, and let $h_p,u_p$ be as in
\eqref{eq:pair}. Then $h_p\in W^{1,2}_{{\textrm{loc}}}(\Omega)\cap C(\Omega)$, $u_p\in C^1(\Omega)$, and,
almost everywhere in $\Omega$,
\begin{equation}
\label{eq:system}
\begin{cases}
\partial_{\bar z}h_p-\mu\,\partial_zh_p-\nu\,\overline{\partial_zh_p}=\alpha h_p+\beta u_p,\\
\partial_zu_p=a\,h_p+b\,u_p ,
\end{cases}
\end{equation}
where
\begin{equation}
\label{eq:ab-alphabeta-bound}
\begin{gathered}
|\mu|+|\nu|\leq k=\frac{p-2}p,\qquad |a|+|b|\leq2(p-1)\qquad\mbox{in }\Omega ,\\
|\alpha|+|\beta|\leq C_E\qquad\mbox{on every }E\subset\subset\Omega ,
\end{gathered}
\end{equation}
where $C_E$ is the constant of \eqref{eq:f-growth} relative to $M=\max\{1,\max_{\overline E}(|u|+|\nabla u|)\}$.
\end{prop}

\begin{proof}
As $h_p=2^{1-p}(F_1-iF_2)$, the regularity of $h_p$ and $u_p$ follows directly from
Lemma \ref{lem:stress}, since $u\in C^1$ and $p\geq2$. Moreover, recalling that
$\partial_zu_p=(p-1)|u|^{p-2}h$, by Young's inequality we obtain
\begin{equation}
\label{eq:up-bound}
|\partial_zu_p|\leq(p-2)|u_p|+|h_p| ,
\end{equation}
and by \eqref{eq:f-growth}\footnote{The bounds \eqref{eq:up-bound} and \eqref{eq:g-bound} take the
place of \cite[(15)]{M}. We point out that their linear dependence on the pair $(h_p,u_p)$ makes the
coefficients $a,b,\alpha,\beta$ bounded, as required by the representation theorem of Bers and Nirenberg.}
\begin{equation}
\label{eq:g-bound}
|g|\ \leq\ C_E\big(|u|^{p-1}+|\nabla u|^{p-1}\big)\ =\ C_E\big(|u_p|+2^{p-1}|h_p|\big)\qquad\mbox{on }E .
\end{equation}
Since $\partial_zu_p$ vanishes on $Z$, and $g$ vanishes almost everywhere on $Z$ by \eqref{eq:on-Z},
by \eqref{eq:ab-alphabeta} we obtain
\begin{equation*}
\partial_zu_p=a\,h_p+b\,u_p\quad\mbox{in }\Omega,\qquad -2^{-p}g=\alpha h_p+\beta u_p\quad\mbox{a.e.\ in }\Omega ,
\end{equation*}
and, since $(|h_p|+|u_p|)^2\leq2(|h_p|^2+u_p^2)$, on $E$
\begin{equation*}
\begin{split}
&|a|+|b|=\frac{|\partial_zu_p|\,(|h_p|+|u_p|)}{|h_p|^2+u_p^2}\leq\frac{(p-1)(|h_p|+|u_p|)^2}{|h_p|^2+u_p^2}\leq2(p-1),\\
&|\alpha|+|\beta|=\frac{|g|\,(|h_p|+|u_p|)}{2^p(|h_p|^2+u_p^2)}\leq\frac{2^{p-1}C_E(|h_p|+|u_p|)^2}{2^p(|h_p|^2+u_p^2)}\leq C_E .
\end{split}
\end{equation*}
The bound for $\mu$ and $\nu$ follows from \eqref{eq:munu}.

Finally, let us derive the equation for $h_p$. Since $h_p=2^{1-p}(F_1-iF_2)$ and
\begin{equation*}
\partial_{\bar z}(F_1-iF_2)=\frac12\big(\div F-i\,\mathrm{curl}\,F\big) ,
\end{equation*}
the equation $\div F=-g$ reads
\begin{equation*}
\mathrm{Re}\,\partial_{\bar z}h_p=-2^{-p}g\qquad\mbox{a.e.\ in }\Omega .
\end{equation*}
In the same way, since $h=\frac12(u_x-iu_y)$ and $\nabla u$ has zero curl, $\partial_{\bar z}h$ is real.
Away from $Z$, where $h_p\neq0$, we have $h=|h_p|^{p'-2}h_p$ with $p'=\frac p{p-1}$, a smooth function of
$h_p$, so that $h\in W^{1,2}_{{\textrm{loc}}}(\Omega\setminus Z)$ and, by the chain rule,
\begin{equation*}
\partial_{\bar z}h=|h_p|^{p'-2}\Big(\frac{p'}2\,\partial_{\bar z}h_p+\frac{p'-2}2\,\frac{h_p}{\overline{h_p}}\,\overline{\partial_zh_p}\Big) .
\end{equation*}
Taking imaginary parts, since $\frac{2-p'}{p'}=\frac{p-2}p$, we have
\begin{equation*}
\mathrm{Im}\,\partial_{\bar z}h_p=\frac{p-2}p\,\mathrm{Im}\Big(\frac{h_p}{\overline{h_p}}\,\overline{\partial_zh_p}\Big).
\end{equation*}
Putting these together, we obtain, almost everywhere away from $Z$,
\begin{equation}
\label{eq:hp-eq}
\partial_{\bar z}h_p=\mu\,\partial_zh_p+\nu\,\overline{\partial_zh_p}+\alpha h_p+\beta u_p .
\end{equation}
The identity \eqref{eq:hp-eq} is true also on $Z$, where both sides are zero almost everywhere
by \eqref{eq:on-Z}.
\end{proof}

In the first equation of \eqref{eq:system} the term $\beta u_p$ does not involve the unknown $h_p$,
while the representation theorem of Bers and Nirenberg requires a homogeneous system: this is the role
of the modified stress field $\mathcal R$, as anticipated in the introduction.

\subsection*{The modified stress field}

We set
\begin{equation*}
\mathcal L=\partial_{\bar z}-\mu\,\partial_z-\nu\,\overline{\partial_z} ,
\end{equation*}
so that, since $u_p$ is real, the system \eqref{eq:system} can be written as
\begin{equation*}
\mathcal Lh_p=\alpha h_p+\beta u_p,\qquad \partial_{\bar z}u_p=\bar a\,\overline{h_p}+\bar b\,u_p .
\end{equation*}
For a disc $B_R=B_R(z_0)\subset\subset\Omega$ and $v\in W^{1,2}(B_R)\cap L^\infty(B_R)$, we consider the
modified stress field
\begin{equation*}
\mathcal R=h_p-v\,u_p .
\end{equation*}
By the product rule and the system, we have
\begin{equation*}
\begin{split}
\mathcal L(v\,u_p)&=u_p\,\mathcal Lv+v\,\partial_{\bar z}u_p-\mu v\,\partial_zu_p-\nu\bar v\,\overline{\partial_zu_p}\\
&=u_p\,\mathcal Lv+\bar a(v-\nu\bar v)\,\overline{h_p}-\mu av\,h_p+\big((\bar b-\mu b)v-\nu\bar b\,\bar v\big)u_p ,
\end{split}
\end{equation*}
and subtracting this from $\mathcal Lh_p$, with $h_p=\mathcal R+v\,u_p$, we obtain, almost everywhere in
$B_R$,
\begin{equation}
\label{eq:riccati}
\mathcal L\mathcal R=(\alpha+\mu av)\,\mathcal R+\bar a(\nu\bar v-v)\,\overline{\mathcal R}
+\big(\mathcal Q(v)-\mathcal Lv\big)u_p ,
\end{equation}
where
\begin{equation}
\label{eq:Q}
\mathcal Q(v)=\beta+\alpha v+(\mu b-\bar b)v+\nu\bar b\,\bar v+\mu a\,v^2+\nu\bar a\,\bar v^2-\bar a|v|^2 .
\end{equation}
Hence $\mathcal R$ solves a homogeneous Beltrami-type system as soon as $v$ solves the Riccati-type
equation
\begin{equation*}
\mathcal Lv=\mathcal Q(v) ,
\end{equation*}
quadratic in $v$, which we may now solve on small discs by a standard Banach--Caccioppoli contraction
mapping argument.

\begin{lemma}[Reduction to a homogeneous system for $\mathcal R$]
\label{lem:remove-W}
Let $u$ be a solution of \eqref{eq:eq-f} and $z_0\in\Omega$. There is $R_0>0$ such that, for
$0<R\leq R_0$ and $B_R=B_R(z_0)$, the modified stress field $\mathcal R=h_p-v\,u_p$ satisfies
\begin{equation}
\label{eq:R-hom}
\partial_{\bar z}\mathcal R=\mu\,\partial_z\mathcal R+\nu\,\overline{\partial_z\mathcal R}+A\mathcal R+B\overline{\mathcal R}
\qquad\mbox{a.e.\ in }B_R ,
\end{equation}
for some $v\in W^{1,q}(B_R)\cap C(\overline{B_R})$, $q>2$, with $|v|\leq1$, where
\begin{equation*}
A=\alpha+\mu av,\qquad B=\bar a(\nu\bar v-v)
\end{equation*}
are bounded on $B_R$. Moreover, $\mathcal R\in W^{1,2}(B_R)\cap C(\overline{B_R})$ and, if $u(z_0)=0$ and
$z_0$ is not an interior point of $\{u=0\}$, $\mathcal R$ does not vanish identically on any disc
$B_r(z_0)$, $0<r\leq R$.
\end{lemma}

\begin{proof}
We argue as in the proof of \cite[Lemma 3.3]{M}. Let $\overline{B_{R_0}}\subset\Omega$ and
$E=B_{R_0}$. By \eqref{eq:riccati}, it suffices to find $v$ with
\begin{equation*}
\mathcal Lv=\mathcal Q(v)\quad\mbox{in }B_R,\qquad |v|\leq1 ,
\end{equation*}
and then $A$ and $B$ are bounded by \eqref{eq:ab-alphabeta-bound}. On the unit ball $X$ of
$L^\infty(B_R)$, by \eqref{eq:ab-alphabeta-bound} and \eqref{eq:Q},
\begin{equation*}
\|\mathcal Q(v)\|_\infty\leq M_{\mathcal Q},\qquad
\|\mathcal Q(v)-\mathcal Q(w)\|_\infty\leq L_{\mathcal Q}\,\|v-w\|_\infty\qquad(v,w\in X),
\end{equation*}
with $M_{\mathcal Q},L_{\mathcal Q}$ depending only on $p$ and $C_E$, while the right inverse of
$\mathcal L$ satisfies, for some $q>2$,
\begin{equation*}
\|\mathcal L^{-1}\phi\|_\infty\leq C_0R\,\|\phi\|_\infty,\qquad
\mathcal L^{-1}\phi\in W^{1,q}(B_R)\cap C(\overline{B_R}) .
\end{equation*}
Hence, for $C_0R\max\{M_{\mathcal Q},L_{\mathcal Q}\}<1$, the map $v\mapsto\mathcal L^{-1}\mathcal Q(v)$ is a
contraction of $X$, and its fixed point is the desired $v$. The regularity of $\mathcal R$ follows from
that of $h_p$ and $v\,u_p$. Finally, if $\mathcal R\equiv0$ on $B_r(z_0)$, then
\begin{equation*}
|h|^{p-1}=|h_p|=|v|\,|u_p|\leq|u|^{p-1},\qquad\mbox{that is}\qquad |\nabla u|\leq2|u|\quad\mbox{on }B_r(z_0).
\end{equation*}
As in \cite[Lemma 3.3]{M}, since $u(z_0)=0$, the fundamental theorem of calculus on
the rays emanating from $z_0$ yields $u\equiv0$ on $B_r(z_0)$, and this concludes the proof.
\end{proof}

By the representation theorem of Bers and Nirenberg, see Bojarski \cite[Theorem 4.4]{Bojarski},
every solution of \eqref{eq:R-hom}, with $|\mu|+|\nu|\leq k<1$ and $A,B$ bounded, is of the form
\begin{equation*}
\mathcal R=\mathcal F\circ\chi ,
\end{equation*}
where $\mathcal F$ is a holomorphic function times a continuous nonvanishing one, and $\chi$ is a
quasiconformal homeomorphism with distortion constant
\begin{equation*}
K=\frac{1+k}{1-k}=p-1 .
\end{equation*}
Hence the zeros of $\mathcal R$ are isolated, with positive integer multiplicities, unless $\mathcal R$
vanishes identically, which Lemma \ref{lem:remove-W} excludes near $z_0$ when $u(z_0)=0$ and $z_0$ is
not an interior point of $\{u=0\}$. These observations, as
in \cite[Lemma 4.1]{M}, yield the following

\begin{lemma}[Representation of $\mathcal R$]
\label{lem:repr-W}
Let $z_0$, $B_R=B_R(z_0)$, $v$ and $\mathcal R=h_p-v\,u_p$ be as in Lemma \ref{lem:remove-W}, and let
$K=p-1$. Then there are a $K$-quasiconformal homeomorphism $\chi$ of $\mathbb C$ onto $\mathbb C$ and
a continuous function $\mathcal F$ on $\chi(B_R)$ such that
\begin{equation*}
\mathcal R=\mathcal F\circ\chi\qquad\mbox{in }B_R .
\end{equation*}
If moreover $u(z_0)=0$ and $z_0$ is not an interior point of $\{u=0\}$, then near each zero $\zeta_1$
of $\mathcal F$ it holds that
\begin{equation*}
\mathcal F(\zeta)=(\zeta-\zeta_1)^n\,\Theta(\zeta),
\end{equation*}
for some integer $n\geq1$ and a nonvanishing continuous function $\Theta$; in particular, the zeros of
$\mathcal R$ in $B_R$ are isolated.
\end{lemma}

\section{Unique continuation in the plane}
\label{sec:uc}

In this section $N=2$, $p>1$, and $u$ is a solution of \eqref{eq:eq-f}; we prove
Theorem \ref{thm:plane}(i)--(ii).

\begin{proof}[Proof of Theorem \ref{thm:plane}(i)--(ii)]
Let $E\subset\subset\Omega$. For every $p>1$, testing the divergence form of \eqref{eq:eq-f} with
$\zeta^pu$, as for $p$-harmonic functions in Lindqvist \cite[Lemma 2.9]{Lindqvist}, and absorbing the
terms due to $f$ by \eqref{eq:f-growth} and Young's inequality, we obtain the Caccioppoli-type inequality: if
$\overline{B_{2\rho}(z_0)}\subset E$ and $2\rho\leq1$, then
\begin{equation}
\label{eq:caccioppoli}
\int_{B_\rho(z_0)}|\nabla u|^p\,dx\ \leq\ \frac C{\rho^p}\int_{B_{2\rho}(z_0)}|u|^p\,dx ,\qquad C=C(p,C_E) .
\end{equation}

Let first $1<p<2$. By \eqref{eq:f-growth} and Young's inequality, on $E$
\begin{equation*}
|f(x,u,\nabla u)|\,|\nabla u|^{2-p}\leq C_E\big(|u|^{p-1}|\nabla u|^{2-p}+|\nabla u|\big)\leq2C_E\big(|u|+|\nabla u|\big) ,
\end{equation*}
which extends to \eqref{eq:eq-f} the bound used in the proof of \cite[Theorem A]{M}. Moreover, the
second derivatives of $u$, needed there, are square summable and vanish almost everywhere on $Z$, by
Lemma \ref{lem:stress}: with $p'=\frac p{p-1}>2$ we have $\nabla u=|F|^{p'-2}F$, so that, by the chain
rule,
\begin{equation*}
\partial_{ij}u=|F|^{p'-2}\big(\delta_{ik}+(p'-2)\,e_ie_k\big)\,\partial_jF_k\quad\mbox{on }\Omega\setminus Z,\qquad
D^2u=0\quad\mbox{a.e.\ on }Z ,
\end{equation*}
summing over $k$, where $e=\frac F{|F|}$. Hence, as anticipated in \cite[Remark 1.1]{M},
the same proof carries over with no modifications whatsoever.

Let now $p\geq2$: the new case.

(i) Let $\Omega_0$ be the interior of the nodal set of $u$, and suppose that
$\emptyset\neq\Omega_0\neq\Omega$. Since $\Omega$ is connected, some $z_0\in\Omega\setminus\Omega_0$ is a
limit of points of $\Omega_0$; then $u(z_0)=0$, and $z_0$ is not an interior point of $\{u=0\}$. By
Lemma \ref{lem:repr-W}, the zeros of $\mathcal R$ in $B_R(z_0)$ are isolated, whereas $\mathcal R$
vanishes on the nonempty open set $\Omega_0\cap B_R(z_0)$, where $h_p=u_p=0$. Hence, if $u$ vanishes on a
nonempty open set, then $\Omega_0=\Omega$.

(ii) Below, $c$ denotes positive constants independent of $\rho$, and we prove \eqref{eq:finite-order}
for $\rho$ small. If $u(z_0)\neq0$, then $|u|\geq\frac12|u(z_0)|$ on $B_\rho(z_0)$, so that
\begin{equation*}
\int_{B_\rho(z_0)}|u|^pdx\geq c\,\rho^2 ,
\end{equation*}
that is $m=2$. If $u(z_0)=0\neq\nabla u(z_0)$, let
$e=\nabla u(z_0)/|\nabla u(z_0)|$ and $D_\rho=B_{\rho/4}(z_0+\frac\rho2\,e)\subset B_\rho(z_0)$. For
$z\in D_\rho$ we have $(z-z_0)\cdot e\geq\frac\rho4$, so that, by Taylor's formula,
\begin{equation*}
u(z)=|\nabla u(z_0)|\,(z-z_0)\cdot e+o(\rho)\geq\frac\rho4\,|\nabla u(z_0)|+o(\rho)\geq\frac\rho8\,|\nabla u(z_0)| ,
\end{equation*}
and therefore
\begin{equation*}
\int_{B_\rho(z_0)}|u|^pdx\geq\int_{D_\rho}|u|^pdx\geq\frac{\pi\rho^2}{16}\Big(\frac\rho8\,|\nabla u(z_0)|\Big)^p=c\,\rho^{p+2} ,
\end{equation*}
that is $m=p+2$. Let finally $z_0\in\Sigma_0$, so that $\mathcal R(z_0)=0$. By (i), $u\not\equiv0$ near
$z_0$, so that, by Lemmas \ref{lem:remove-W} and \ref{lem:repr-W}, near $z_0$
\begin{equation*}
\mathcal R=\mathcal F\circ\chi,\qquad \mathcal F(\zeta)=(\zeta-\zeta_0)^n\,\Theta(\zeta),\qquad \zeta_0=\chi(z_0),
\end{equation*}
with $|\Theta|\geq c$. By the H\"older distortion of quasiconformal mappings \cite[Lemma 4.2]{M},
based on Mori's theorem, $|\chi(z)-\zeta_0|\geq c\,|z-z_0|^K$ with $K=p-1$, and therefore
\begin{equation*}
|\mathcal R(z)|=|\chi(z)-\zeta_0|^n\,|\Theta(\chi(z))|\geq c\,|z-z_0|^{nK}\qquad\mbox{on }B_{\rho_0}(z_0) ,
\end{equation*}
where we may assume $\overline{B_{2\rho_0}(z_0)}\subset\Omega$ and $2\rho_0\leq1$. Since $Kp'=p$,
integrating we obtain, for $0<\rho<\rho_0$,
\begin{equation*}
\int_{B_\rho(z_0)}|\mathcal R|^{p'}dx\geq c\int_{B_\rho(z_0)}|z-z_0|^{np}\,dx=\frac{2\pi c}{np+2}\,\rho^{np+2} .
\end{equation*}
On the other hand, since $|v|\leq1$ and $p'\leq2\leq p$, by \eqref{eq:pair} we have
\begin{equation*}
|\mathcal R|^{p'}\leq2^{p'-1}\big(|h_p|^{p'}+|u_p|^{p'}\big)=2^{p'-1}\big(2^{-p}|\nabla u|^p+|u|^p\big)
\leq|\nabla u|^p+2|u|^p .
\end{equation*}
Hence, by \eqref{eq:caccioppoli} with $E=B_{2\rho_0}(z_0)$, since $\rho\leq1$,
\begin{equation*}
c\,\rho^{np+2}\leq\int_{B_\rho(z_0)}\big(|\nabla u|^p+2|u|^p\big)dx\leq\frac{C+2}{\rho^p}\int_{B_{2\rho}(z_0)}|u|^p\,dx ,
\end{equation*}
which is \eqref{eq:finite-order} with $m=np+p+2$, up to renaming $2\rho$ as $\rho$. Part (iii) is
proved in Section \ref{sec:index}.
\end{proof}

\section{Isolated simultaneous zeros and the index of the gradient}
\label{sec:index}

In this section we prove Theorem \ref{thm:plane}(iii). Throughout, $N=2$, $p>1$, $\Omega$ is
connected, $u\not\equiv0$ is a solution of \eqref{eq:eq-f}, $z_0\in\Sigma_0$, and the index and the
circles $\gamma_r$ are as in Section \ref{sec:intro}.

We treat both ranges of $p$ at once, setting
\begin{equation*}
(\tilde h,\tilde u)=(h,u)\quad(1<p<2),\qquad (\tilde h,\tilde u)=(h_p,u_p)\quad(p\geq2),
\end{equation*}
and
\begin{equation*}
\mathcal R=\tilde h-v\,\tilde u ,
\end{equation*}
with $v$ as in \cite[Lemma 3.3]{M}, respectively Lemma \ref{lem:remove-W}. In both cases $|v|\leq1$,
$\mathcal R(z_0)=0$ and, by Theorem \ref{thm:plane}(i) and the representation of \cite[Lemma 4.1]{M},
respectively Lemma \ref{lem:repr-W}, on a disc $B_R(z_0)$
\begin{equation*}
\mathcal R=\mathcal F\circ\chi,\qquad \mathcal F(\zeta)=(\zeta-\zeta_0)^n\,\Theta(\zeta),
\end{equation*}
where $\chi$ is a $K$-quasiconformal homeomorphism of $\mathbb C$ onto $\mathbb C$, namely, as observed
in \cite[Section 4]{M}, a homeomorphism which, at small scales, enlarges distances in one direction at most $K$ times more than in
another, with $K=\max\{p-1,\frac1{p-1}\}$ as for $p$-harmonic functions, see Iwaniec and
Manfredi \cite[(13)]{IwaniecManfredi}, $\zeta_0=\chi(z_0)$, $n\geq1$ is an integer, and
$c_1\leq|\Theta|\leq c_2$ on a disc $B_\sigma(\zeta_0)\subset\chi(B_R(z_0))$, for some constants
$0<c_1\leq c_2$. Since $\tilde h$ is a positive multiple of $h$ away from its zeros, and $\nabla u$
corresponds to $2\bar h$, it is enough to show that $\tilde h$ does not vanish on the circles
$\partial B_r(z_0)$, for $r$ small, and that
\begin{equation*}
\mathrm{Ind}_{\tilde h\circ\gamma_r}(0)=n .
\end{equation*}

For $p$-harmonic functions, $h$ is a power of a quasiconformal map near a critical point by Stoilow's
factorisation theorem, and its winding number readily follows, see \cite[(15)]{IwaniecManfredi}. Here only $\mathcal R$ has this
form, and $\tilde h$ differs from it by $v\,\tilde u$. However, a quasiconformal map is
controlled in two respects: as seen in the proof of Theorem \ref{thm:plane}(ii), it
contracts distances no faster than a power and, again by Mori, it distorts small circles only boundedly.
Hence, on each small circle, $|\mathcal R|$ is comparable with its maximum on the disc, while $\tilde u$
is a small fraction of this maximum, see Lemma \ref{lem:u-small} below, so that
\begin{equation*}
|v\,\tilde u|<|\mathcal R|\qquad\mbox{on }\partial B_r(z_0),
\end{equation*}
and $\tilde h=\mathcal R+v\,\tilde u$ has the same winding number as $\mathcal R$.

For $r>0$ such that
\begin{equation}
\label{eq:r-small}
\overline{B_r(z_0)}\subset\chi^{-1}\big(B_\sigma(\zeta_0)\big),
\end{equation}
set
\begin{equation*}
S(r)=\max_{\overline{B_r(z_0)}}|\mathcal R|,\qquad m(r)=\min_{\partial B_r(z_0)}|\mathcal R| .
\end{equation*}
By Mori's theorem on circular distortion, see Astala, Iwaniec and
Martin \cite[(3.45)]{AstalaIwaniecMartin}, $\chi$ maps $\partial B_r(z_0)$ into an annulus centred at
$\zeta_0$ whose radii have ratio at most $C_K$, a constant depending only on $K$. Since $\chi$ is open,
$|\chi-\zeta_0|$ attains its maximum on $\overline{B_r(z_0)}$ at the boundary, and therefore
\begin{equation*}
|\chi(z)-\zeta_0|\leq C_K\,|\chi(z')-\zeta_0|\qquad\mbox{for } z\in\overline{B_r(z_0)},\ z'\in\partial B_r(z_0).
\end{equation*}
By the bounds on $\Theta$, it follows that
\begin{equation*}
|\mathcal R(z)|\leq c_2\,|\chi(z)-\zeta_0|^n\leq c_2\,C_K^n\,|\chi(z')-\zeta_0|^n
\leq\frac{c_2}{c_1}\,C_K^n\,|\mathcal R(z')| ,
\end{equation*}
hence
\begin{equation}
\label{eq:S-m}
S(r)\leq C_*\,m(r),\qquad C_*=\frac{c_2}{c_1}\,C_K^n .
\end{equation}

Together with \eqref{eq:S-m}, the following bound, based on a classical Gronwall-type argument, says
that $v\,\tilde u$ is negligible compared to $\mathcal R$ on small circles.

\begin{lemma}[Smallness of $\tilde u$]
\label{lem:u-small}
For every $r>0$ satisfying \eqref{eq:r-small}, it holds that
\begin{equation}
\label{eq:u-small}
\max_{\overline{B_r(z_0)}}|\tilde u|\leq\eps(r)\,S(r),\qquad\mbox{where}\quad
\eps(r)=\big(2re^{2r}\big)^\tau\quad\mbox{and}\quad\tau=\max\{1,p-1\} .
\end{equation}
\end{lemma}

\begin{proof}
Since $|\tilde h|=|h|^\tau$, $|\tilde u|=|u|^\tau$, $|v|\leq1$ and $\tau\geq1$, on $B_r(z_0)$ we have
\begin{equation*}
|\nabla u|=2|h|\leq2\big(S(r)+|u|^\tau\big)^{\frac1\tau}\leq2S(r)^{\frac1\tau}+2|u| .
\end{equation*}
On the rays emanating from $z_0$, the function $\varphi(t)=|u(z_0+t\theta)|$, $|\theta|=1$, vanishes
at $t=0$ and, by the fundamental theorem of calculus, satisfies
\begin{equation*}
\varphi(t)\leq\int_0^t|\nabla u(z_0+s\theta)|\,ds\leq2S(r)^{\frac1\tau}\,t+2\int_0^t\varphi(s)\,ds ,
\end{equation*}
so that, by Gronwall's inequality,
\begin{equation*}
\varphi(t)\leq2t\,e^{2t}\,S(r)^{\frac1\tau} .
\end{equation*}
Hence $|u|\leq2re^{2r}S(r)^{\frac1\tau}$ on $\overline{B_r(z_0)}$, and \eqref{eq:u-small} follows by
raising to the power $\tau$.
\end{proof}

\begin{proof}[Proof of Theorem \ref{thm:plane}(iii)]
Let $r>0$ satisfy \eqref{eq:r-small}. On $\partial B_r(z_0)$, by Lemma \ref{lem:u-small}
and \eqref{eq:S-m},
\begin{equation*}
|\tilde h-\mathcal R|=|v|\,|\tilde u|\ \leq\ \eps(r)\,S(r)\ \leq\ \eps(r)\,C_*\,|\mathcal R| .
\end{equation*}
For $r$ small, $\eps(r)C_*<1$, so that $|\tilde h-\mathcal R|<|\mathcal R|$ on $\partial B_r(z_0)$:
$\tilde h$ does not vanish there, and neither does the homotopy $(1-t)\mathcal R+t\tilde h$, $0\leq t\leq1$.
By the homotopy invariance of the index, $\tilde h\circ\gamma_r$ and $\mathcal R\circ\gamma_r$ have the same
index about $0$. Moreover, $\mathcal R=(\chi-\zeta_0)^n\,\Theta\circ\chi$, where $\chi\circ\gamma_r$ winds
once around $\zeta_0$, since $\chi$ is an orientation-preserving homeomorphism, while
$\Theta\circ\chi\circ\gamma_r$ does not wind around $0$, since $\Theta\circ\chi$ is continuous and
nonvanishing on the disc $\overline{B_r(z_0)}$. Hence, as the index is additive under products, we
obtain
\begin{equation*}
\mathrm{Ind}_{\tilde h\circ\gamma_r}(0)=\mathrm{Ind}_{\mathcal R\circ\gamma_r}(0)
=n\,\mathrm{Ind}_{\chi\circ\gamma_r}(\zeta_0)+\mathrm{Ind}_{\Theta\circ\chi\circ\gamma_r}(0)=n .
\end{equation*}
For $r_0$ small these circles cover $B_{r_0}(z_0)\setminus\{z_0\}$, so that $h$ does not vanish there,
and $z_0$ is an isolated point of $Z$. Finally, on $\partial B_r(z_0)$ the field $\tilde h$ is a positive
multiple of $h$, while $\nabla u$, identified with $u_x+iu_y$, equals $2\bar h$, and complex conjugation
reverses the orientation of closed curves. Therefore
\begin{equation*}
\mathrm{ind}_{z_0}(\nabla u)=\mathrm{Ind}_{\bar h\circ\gamma_r}(0)=-\mathrm{Ind}_{\tilde h\circ\gamma_r}(0)=-n .
\end{equation*}
\end{proof}

\section{Axial symmetry}
\label{sec:axial}

The rest of the paper is devoted to the proof of the main results under axial symmetry, and for this
reason we assume $N\geq3$. In this section we prove Theorem \ref{thm:axial}, providing the details of an
observation in \cite[Section 1.4]{M} on a reduction to lower dimension. In fact, we prove that the
profile of an axially symmetric solution, $u(x)=U(r,x_N)$, solves in the meridian half-plane $\{r>0\}$ a
weighted equation which is equivalent to a two-dimensional $p$-Laplace equation with a first-order drift
term. Away from the axis this term is controlled by $|\nabla U|^{p-1}$, a context where
Theorem \ref{thm:plane} directly applies.

For an axially symmetric domain $\Omega\subset\R^N$, let
\begin{equation*}
\mathcal M=\big\{(r,x_N):\ r>0,\ (r,0,\dots,0,x_N)\in\Omega\big\}
\end{equation*}
be the section of $\Omega\setminus A$ by a meridian half-plane, on which
$\nabla U=(\partial_rU,\partial_{x_N}U)$ denotes the gradient of the profile. Rotating $\mathcal M$ about
the axis, we generate $\Omega\setminus A$: its sections by the meridian half-planes are the leaves of a
foliation, to which the quasiconformal reduction principle can be applied. This is the content of the
following

\begin{prop}[The equation on a leaf]
\label{prop:reduction}
Under the assumptions of Theorem \ref{thm:axial}, the set $\mathcal M$ is connected,
$U\in C^{1,\alpha}_{{\textrm{loc}}}(\mathcal M)$, and
\begin{equation*}
u(x)=U(r,x_N),\qquad|\nabla u(x)|=|\nabla U(r,x_N)|\qquad\mbox{for }x\in\Omega\setminus A .
\end{equation*}
Moreover, there is a measurable function $\tilde f$ on $\mathcal M$ such that
\begin{equation*}
f(x,u,\nabla u)=\tilde f(|x'|,x_N)\qquad\mbox{for a.e.\ }x\in\Omega ,
\end{equation*}
and $U$ is a weak solution in $\mathcal M$ of the weighted equation
\begin{equation}
\label{eq:reduced-div}
-\div\big(r^{N-2}|\nabla U|^{p-2}\nabla U\big)=r^{N-2}\tilde f ,
\end{equation}
and hence of the planar $p$-Laplace equation with a drift term
\begin{equation}
\label{eq:reduced}
-\Delta_pU=\tilde f+\frac{N-2}r\,|\nabla U|^{p-2}\partial_rU ,
\end{equation}
where $\div$ and $\Delta_p$ act on $(r,x_N)$. The right-hand side of \eqref{eq:reduced} satisfies
\eqref{eq:g-bound} on every $E\subset\subset\mathcal M$, with $U$ in place of $u$.
\end{prop}

\begin{proof}
As $A$ has codimension $N-1\geq2$, the set $\Omega\setminus A$ is connected, and so is its image $\mathcal M$
under the continuous map $x\mapsto(|x'|,x_N)$. The profile
\begin{equation*}
U(r,x_N)=u(r,0,\dots,0,x_N)
\end{equation*}
belongs to $C^{1,\alpha}_{{\textrm{loc}}}(\mathcal M)$.

Away from $A$ we write $e_i=\frac{x_i}r$, $i=1,\dots,N-1$, so that
\begin{equation*}
\sum_{i=1}^{N-1}e_i^2=1,\qquad \partial_ir=\frac{x_i}r=e_i,\qquad
\partial_ie_i=\partial_i\Big(\frac{x_i}r\Big)=\frac1r-\frac{x_i}{r^2}\,\partial_ir=\frac1r-\frac{x_i^2}{r^3}=\frac{1-e_i^2}r .
\end{equation*}
We have
\begin{equation*}
\partial_iu=\partial_rU\,\partial_ir=\partial_rU\,e_i\quad(i\leq N-1),\qquad \partial_Nu=\partial_{x_N}U ,
\end{equation*}
that is
\begin{equation*}
\nabla u=(\partial_rU\,e_1,\dots,\partial_rU\,e_{N-1},\partial_{x_N}U),\qquad
|\nabla u|^2=(\partial_rU)^2\sum_{i=1}^{N-1}e_i^2+(\partial_{x_N}U)^2=|\nabla U|^2 ,
\end{equation*}
and the stress field of $u$ is
\begin{equation*}
F=|\nabla u|^{p-2}\nabla u=(a\,e_1,\dots,a\,e_{N-1},\,b),\qquad
a=|\nabla U|^{p-2}\partial_rU,\quad b=|\nabla U|^{p-2}\partial_{x_N}U .
\end{equation*}
By Lemma \ref{lem:stress}, $F\in W^{1,2}_{{\textrm{loc}}}(\Omega)$, so that $a=\sum_{i=1}^{N-1}F_ie_i$ and $b=F_N$
belong to $W^{1,2}_{{\textrm{loc}}}$ away from $A$, and, as functions of $(r,x_N)$, to
$W^{1,2}_{{\textrm{loc}}}(\mathcal M)$, since $r^{N-2}$ is bounded and bounded away from zero on compact
subsets of $\mathcal M$. With this regularity at hand, and since $a$ depends on $x'$ only through $r$, so
that $\partial_ia=\partial_ra\,e_i$, away from $A$ we obtain
\begin{equation*}
\begin{split}
\sum_{i=1}^{N-1}\partial_i(a\,e_i)
&=\sum_{i=1}^{N-1}\big(\partial_ia\;e_i+a\,\partial_ie_i\big)
=\partial_ra\sum_{i=1}^{N-1}e_i^2+\frac ar\sum_{i=1}^{N-1}\big(1-e_i^2\big)\\
&=\partial_ra+\frac{N-2}r\,a ,
\end{split}
\end{equation*}
hence
\begin{equation*}
\div F=\partial_ra+\partial_{x_N}b+\frac{N-2}r\,a .
\end{equation*}
As $-\div F=g$ almost everywhere by Lemma \ref{lem:stress}, the function $g=f(x,u,\nabla u)$ coincides
almost everywhere with a function of $(|x'|,x_N)$, which we write as $\tilde f(|x'|,x_N)$. Since moreover
\begin{equation*}
\partial_r\big(r^{N-2}a\big)+\partial_{x_N}\big(r^{N-2}b\big)
=r^{N-2}\Big(\partial_ra+\partial_{x_N}b+\frac{N-2}r\,a\Big),
\end{equation*}
the equation $-\div F=g$ reads, in the variables $(r,x_N)$, as \eqref{eq:reduced-div}, or equivalently
as \eqref{eq:reduced}, almost everywhere and hence, as $a$ and $b$ belong to
$W^{1,2}_{{\textrm{loc}}}(\mathcal M)$, in the weak sense. Finally, on a set $E$ with $\overline E\subset \mathcal M$ compact,
where $r\geq r_E>0$, by \eqref{eq:f-growth}
\begin{equation*}
\Big|\tilde f+\frac{N-2}r\,|\nabla U|^{p-2}\partial_rU\Big|\ \leq\ C_E\big(|U|^{p-1}+|\nabla U|^{p-1}\big)+\frac{N-2}{r_E}\,|\nabla U|^{p-1},
\end{equation*}
as $|u|=|U|$ and $|\nabla u|=|\nabla U|$.
\end{proof}

\begin{proof}[Proof of Theorem \ref{thm:axial}]
By the above proposition Theorem \ref{thm:axial} follows from Theorem \ref{thm:plane} applied to the
profile $U$ on the connected set $\mathcal M$, as the proof of Theorem \ref{thm:plane} uses
\eqref{eq:f-growth} only through the bound \eqref{eq:g-bound} on compact subsets. Observe that the critical points and the simultaneous zeros of $u$
away from the axis are now the orbits of those of $U$. As $x\mapsto(|x'|,x_N)$ is an open map on
$\Omega\setminus A$, (i) and (iii) of Theorem \ref{thm:axial} follow then from (i) and (iii) of
Theorem \ref{thm:plane} for $U$, by definition of isolated orbits in the case of (iii).

As for (ii), let $r_0=|x_0'|$. Integrating over $B_\rho(Gx_0)$, which corresponds to the disc
$B_\rho(r_0,x_{0,N})$, in cylindrical coordinates, where $dx=r^{N-2}\,dr\,d\sigma(\omega)\,dx_N$ and
$r\geq\frac{r_0}2$ for $\rho$ small, we obtain
\begin{equation*}
\int_{B_\rho(Gx_0)}|u|^p\,dx\ \geq\ c\int_{B_\rho(r_0,x_{0,N})}|U|^p\,dr\,dx_N ,\qquad
c=|S^{N-2}|\Big(\frac{r_0}2\Big)^{N-2} ,
\end{equation*}
and (ii) follows from \eqref{eq:finite-order} for $U$, with the same $m$.
\end{proof}

\section{Liouville-type theorems}
\label{sec:nonex}

As in the classical paper of Esteban and
Lions \cite{EstebanLions}, a Pohozaev identity makes the gradient of a solution vanish on a piece of the
boundary, across which the solution extends by zero, and weak unique continuation
forces it to be flat zero everywhere.

The proofs of \cite[Theorems B and C]{M}, including the regularity, apply in every dimension, and unique
continuation is now provided by Theorem \ref{thm:axial}(i), which reaches the degenerate regime $p\geq2$. On the
half-space, $\partial_Nu=0$ on $\partial H$ by the result of Farina, the author and
Willem \cite[Lemma 2.1]{FarinaMercuriWillem}, valid for every $1<p<N$, see also \cite[Lemma 6.1]{M},
\cite[Lemma 2.3]{MercuriWillem} and, for $p=2$, \cite[proof of Theorem I.1]{EstebanLions}; on a domain star-shaped with respect to $x_0$, the Pohozaev
identity of Guedda and V\'eron \cite[Theorem 1.1]{GueddaVeron} yields $\nabla u=0$ on the nonempty,
relatively open set
\begin{equation*}
\Gamma=\{x\in\partial\Omega:\ (x-x_0)\cdot\nu(x)>0\} .
\end{equation*}

\begin{proof}[Proof of Theorem \ref{thm:nonex}]
(i) The extension of $u$ by zero to $\R^N$ is axially symmetric, as $A$ is normal to $\partial H$, and it
vanishes on $\{x_N<0\}$. Since $f(s)=|s|^{p^*-2}s$ satisfies \eqref{eq:f-growth} with $C_M=M^{p^*-p}$,
Theorem \ref{thm:axial}(i) gives $u\equiv0$.

(ii) Let $B$ be a ball centred at a point of $\Gamma$, in which $\partial\Omega$ is a graph contained in $\Gamma$, and
$GB=\{gx:\ g\in G,\ x\in B\}$. By symmetry, $u$ and $\nabla u$ vanish on $\partial\Omega\cap GB$, which is
covered by the balls $gB$, and the extension of $u$ by zero is an axially symmetric solution in
$\Omega\cup GB$, vanishing on the open set $GB\setminus\overline\Omega$. Theorem \ref{thm:axial}(i) gives
$u\equiv0$.
\end{proof}

\section{Blow-up analysis for critical problems}
\label{sec:compact}

In this section we prove Theorem \ref{thm:compact}, which extends the classical result of
Struwe \cite{Struwe} and its version for symmetric domains due to Clapp \cite{Clapp}, in the spirit of
P.-L. Lions \cite{LionsSym}, and refines its quasilinear generalisations by the author and
Willem \cite[Theorem 1.2]{MercuriWillem} and, for symmetric Palais--Smale sequences, by the author and
Pacella \cite[Theorem 2.2]{MercuriPacella}. The method follows the scheme of
Willem \cite[Theorem 8.13]{Willem} for $p=2$ and of \cite[Theorem 1.2]{MercuriWillem} for the
$p$-Laplacian, with the boundedness of the Palais--Smale sequences justified as in \cite[proof of
Theorem D]{M}. We highlight only the new aspects. The bubbles concentrate on the axis, by an
observation of Clapp \cite{Clapp} which here takes a particularly simple form, and the profiles on half-spaces,
which for sign-changing sequences could not be excluded in \cite[Remark 2.5]{MercuriPacella}, are
ruled out by Theorem \ref{thm:nonex}(i).

\begin{proof}[Proof of Theorem \ref{thm:compact}]
The boundedness of the Palais--Smale sequences is proved as in \cite[proof of Theorem D]{M}. The
iteration then proceeds as in Clapp \cite{Clapp} for $p=2$, and as in \cite[Theorem 2.2]{MercuriPacella}
for $p\neq2$, and stops after finitely many steps: along a subsequence, we obtain an axially
symmetric solution $v_0$, closed subgroups $\Gamma_1,\dots,\Gamma_k$ of $G$ with $G/\Gamma_i$ finite,
nontrivial solutions $v_i$ of
\begin{equation*}
-\Delta_pu=|u|^{p^*-2}u ,
\end{equation*}
invariant under $\Gamma_i$ and defined on $\R^N$ or on a half-space $H_i$, scales $\eps_n^i$ tending to
zero, and points $y_n^i\in\Omega$ with isotropy group $\Gamma_i$, such that
\begin{gather*}
\Big\|u_n-v_0-\sum_{i=1}^k\sum_{[g]\in G/\Gamma_i}(\eps_n^i)^{\frac{p-N}p}\,
v_i\Big(\frac{g^{-1}(\cdot-gy_n^i)}{\eps_n^i}\Big)\Big\|\to0 ,\\
\|u_n\|^p\to\|v_0\|^p+\sum_{i=1}^k|G/\Gamma_i|\,\|v_i\|^p,\\
J_\Omega(v_0)+\sum_{i=1}^k|G/\Gamma_i|\,J_{\R^N}(v_i)=c ,
\end{gather*}
where $G/\Gamma_i$ is identified with the $G$-orbit of $y_n^i$, the class $[g]$ corresponding to the point
$gy_n^i$, $|G/\Gamma_i|$ is the cardinality of this orbit, and the $v_i$ are extended by zero outside $H_i$.

\emph{Concentration on the axis.} The $G$-orbit of $y_n^i$, of cardinality $|G/\Gamma_i|$, is finite,
while the $G$-orbit of a point away from the axis is a sphere of positive
dimension $N-2$. Hence $\Gamma_i=G$ and the points $y_n^i$ lie on the axis; each sum over $G/\Gamma_i$
reduces to one term, and the profiles $v_i$ are axially symmetric.

\emph{No half-space profile.} Let $H_i$ be, by contradiction, one of the half-spaces above. This case
occurs when
\begin{equation*}
\sup_n\frac{\operatorname{dist}(y_n^i,\partial\Omega)}{\eps_n^i}<\infty .
\end{equation*}
As the limit of the rescaled domains $(\Omega-y_n^i)/\eps_n^i$, which are axially
symmetric, $H_i$ is axially symmetric, and its boundary is orthogonal to the axis. By a rigid motion
commuting with $G$, we can assume that $H_i=H$, and $v_i\in\DR_0^{1,p}(H)$ is an axially symmetric
solution to
\begin{equation*}
-\Delta_pu=|u|^{p^*-2}u\qquad\mbox{in }H .
\end{equation*}
However, by Theorem \ref{thm:nonex}(i), $v_i\equiv0$. Hence
\begin{equation*}
\frac{\operatorname{dist}(y_n^i,\partial\Omega)}{\eps_n^i}\to\infty .
\end{equation*}
Finally, if $\Omega$ does not meet the axis, there are no points $y_n^i$, so that $k=0$ and $u_n\to v_0$
strongly in $W^{1,p}_0(\Omega)$. This concludes the proof.
\end{proof}

\section*{Acknowledgements}
Carlo Mercuri is a member of the group GNAMPA of Istituto Nazionale di Alta Matematica (INdAM).

\section*{Statements and Declarations}

\noindent\textbf{Funding.} The author did not receive support from any organisation for the submitted work.

\noindent\textbf{Competing interests.} The author has no relevant financial or non-financial interests to disclose.

\end{document}